\documentclass{elsart3-1-ppt}

\usepackage{amssymb,amsmath,latexsym,setspace}
\usepackage[english,francais]{babel}

\newtheorem{e-proposition}[theorem]{Proposition}

\newtheorem{e-definition}[theorem]{Definition\rm}

\renewcommand{\theequation}{\arabic{equation}}
\renewenvironment{thebibliography}[1]{\begin{oldthebibliography}{#1}\setlength{\itemsep}{-.5ex}}{\end{oldthebibliography}}

\newcommand{\M}{\mathfrak M}
\newcommand{\iM}[1]{\int_{\M}{#1}\,dv_g}
\newcommand{\nrM}[2]{\|{#1}\|_{\mathrm L^{#2}(\M)}}

\def\og{\leavevmode\raise.3ex\hbox{$\scriptscriptstyle\langle\!\langle$~}}
\def\fg{\leavevmode\raise.3ex\hbox{~$\!\scriptscriptstyle\,\rangle\!\rangle$}}

\newcommand{\R}{{\mathbb R}}
\renewcommand{\S}{{\mathbb S}}

\newcommand{\be}[1]{\begin{equation}\label{#1}}
\newcommand{\ee}{\end{equation}}
\renewcommand{\(}{\left(}
\renewcommand{\)}{\right)}

\usepackage{color}

\renewcommand{\S}{\mathbb{S}}

\newcommand{\nrmRd}[2]{\|{#1}\|_{\mathrm L^{#2}(\R^d)}}

\newcommand{\finprf}{\unskip\null\hfill$\;\square$\vskip 0.3cm}
\newenvironment{proof}{\par\noindent{\emph{Proof. }}}{\finprf}

\newcommand{\irdmu}[1]{\int_{\R^+\times\M}#1\,d\mu}

\newcommand{\nrmcndgen}[2]{\|{#1}\|_{\mathrm L^{#2}(\mathcal C)}}

\newcommand{\iR}[1]{\int_{\R}{#1}\,ds}
\newcommand{\nrml}[2]{\|{#1}\|_{\mathrm L^{#2}(\R)}}

\newcommand{\alf}{\lambda}

\journal{the Acad\'emie des sciences}
\begin{document}

\centerline{}
\begin{frontmatter}

\selectlanguage{english}

%%%%%%%%%%%%%%%%%%%%%%%%%%%%%%%%%%%%%%%%%%%%%%%%%%%%%%%%%%%%%%%%%%%%%%%%
%%%%%%%%%%%%%%%%%%%%%%%%%%%%%%%%%%%%%%%%%%%%%%%%%%%%%%%%%%%%%%%%%%%%%%%%
\title{Keller-Lieb-Thirring inequalities\\ for Schr\"odinger operators on cylinders}

\selectlanguage{english}
\author[Ceremade]{Jean Dolbeault}
\ead{dolbeaul@ceremade.dauphine.fr}
\and
\author[Ceremade]{Maria J.~Esteban}
\ead{esteban@ceremade.dauphine.fr}
\and
\author[GIT]{Michael Loss}
\ead{loss@math.gatech.edu}

\address[Ceremade]{Ceremade (UMR CNRS no. 7534), Universit\'e Paris-Dauphine, Place de Lattre de Tassigny, 75775 Paris 16, France}
\address[GIT]{Skiles Building, Georgia Institute of Technology, Atlanta GA 30332-0160, USA}

\begin{abstract}\selectlanguage{english}
This note is devoted to Keller-Lieb-Thirring spectral estimates for Schr\"odinger operators on infinite cylinders: the absolute value of the ground state level is bounded by a function of a norm of the potential. Optimal potentials with small norms are shown to depend on a single variable: this is a \emph{symmetry result}. The proof is a perturbation argument based on recent rigidity results for nonlinear elliptic equations on cylinders. Conversely, optimal single variable potentials with large norms must be unstable: this provides a \emph{symmetry breaking result}. The optimal threshold between the two regimes is established in the case of the product of a sphere by a line.
\vskip 0.5\baselineskip\noindent
{\bf Keywords.} Cylinder; Schr\"odinger operator; eigenvalues; Keller-Lieb-Thirring inequalities; optimal constants; symmetry; symmetry breaking
% Sobolev inequality; Gagliardo-Nirenberg inequality; interpolation; rigidity results; fast diffusion equation; Laplace-Beltrami operator; spectral estimates; bifurcation; instability
\\[4pt]
\noindent MSC (2010):
58J50;
% 58J50 Spectral problems; spectral geometry; scattering theory [See also 35Pxx]
81Q10;
% 81Q10 Selfadjoint operator theory in quantum theory, including spectral analysis
81Q35;
% 81Q35 Quantum mechanics on special spaces: manifolds, fractals, graphs, etc.
35P15
\selectlanguage{francais}
\vskip 0.5\baselineskip \noindent {\bf In\'egalit\'es de Keller-Lieb-Thirring pour des op\'erateurs de Schr\"odinger sur des cylindres.}

\vskip 0.5\baselineskip\noindent{\bf R\'esum\'e.}
Cette note est consacr\'ee \`a des estimations spectrales de Keller-Lieb-Thirring pour des op\'erateurs de Schr\"odinger sur des cylindres infinis: la valeur absolue de l'\'etat fondamental est born\'ee par une fonction d'une norme du potentiel. Il est montr\'e que les potentiels optimaux de petite norme ne d\'ependent que d'une seule variable: il s'agit d'un \emph{r\'esultat de sym\'etrie}. La preuve provient d'un argument de perturbation qui repose sur des r\'esultats de rigidit\'e r\'ecents pour des \'equations elliptiques non-lin\'eaires sur des cylindres. A l'inverse, les potentiels optimaux de grande norme qui ne d\'ependent que d'une seule variable sont instables: cela fournit un \emph{r\'esultat de brisure de sym\'etrie}. La valeur optimale qui s\'epare les deux r\'egimes est \'etablie dans le cas du produit d'une sph\`ere et d'une droite.
\end{abstract}
\end{frontmatter}\vspace*{-1.5cm}
%%%%%%%%%%%%%%%%%%%%%%%%%%%%%%%%%%%%%%%%%%%%%%%%%%%%%%%%%%%%%%%%%%%%%%%%
%%%%%%%%%%%%%%%%%%%%%%%%%%%%%%%%%%%%%%%%%%%%%%%%%%%%%%%%%%%%%%%%%%%%%%%%
\selectlanguage{english}
\setcounter{equation}{0}
%%%%%%%%%%%%%%%%%%%%%%%%%%%%%%%%%%%%%%%%%%%%%%%%%%%%%%%%%%%%%%%%%%%%%%%%
%%%%%%%%%%%%%%%%%%%%%%%%%%%%%%%%%%%%%%%%%%%%%%%%%%%%%%%%%%%%%%%%%%%%%%%%

%%%%%%%%%%%%%%%%%%%%%%%%%%%%%%%%%%%%%%%%%%%%%%%%%%%%%%%%%%%%%%%%%%%%%%%%
%%%%%%%%%%%%%%%%%%%%%%%%%%%%%%%%%%%%%%%%%%%%%%%%%%%%%%%%%%%%%%%%%%%%%%%%
\section{Introduction and main results}\label{Sec:Intro}

Let $(\M,g)$ be a smooth compact connected Riemannian manifold of dimension $d-1$, without boundary. We denote by~$ \mathrm{Ric}$ the Ricci tensor and by~$\lambda_1^\M$ the lowest positive eigenvalue of the Laplace-Beltrami operator $-\,\Delta_g$ on $\M$. Let us define the minimum of the Ricci curvature of $\M$ by $\kappa:=\inf_{\M}\inf_{\xi\in\S^{d-2}}\mathrm{Ric}(\xi\,,\xi)$ and consider the infinite cylinder $\mathcal C:=\R\times\M$. We shall denote by $x=(s,z)$ the variable on $\mathcal C$, so that the Laplace-Beltrami operator on $\mathcal C$ can be written as $-\,(\partial^2_s+\Delta_g)$. For simplicity, we shall assume that $\mathrm{vol}_g(\M)=1$, so that $\nrmcndgen Vq=\nrml Vq$ if $V$ is a potential that depends only on $s$. The goal of this note is to compare
\[
\Lambda(\mu):=\sup\left\{\lambda_1^{\mathcal C}[V]:V\in\mathrm L^q(\mathcal C)\,,\;\nrmcndgen Vq=\mu\right\}\;\mbox{and}\;\Lambda_\R(\mu):=\sup\left\{\lambda_1^{\R}[V]:V\in\mathrm L^q(\R)\,,\;\nrml Vq=\mu\right\}
\]
where $-\lambda_1^{\mathcal C}[V]$ and $-\lambda_1^{\R}[V]$ denote the lowest eigenvalues of the Schr\"odinger operators $-\partial^2_s-\,\Delta_g-V$ and $-\partial^2_s-V$ respectively on $\mathcal C$ and $\R$. What we aim at is a \emph{symmetry result} that allows us to characterize the regime in which optimal potentials depend only on $s$. The expression of $\Lambda_\R(\mu)$ was found by J.B.~Keller in~\cite{MR0121101} and later rediscovered by E.H.~Lieb and W.~Thirring in~\cite{Lieb-Thirring76}. We refer to this result as the \emph{Keller-Lieb-Thirring inequality}, and to~\cite{DEL2011} for its use in the context of Caffarelli-Kohn-Nirenberg inequalities. Let us define
\[
\textstyle\mu_1:=q\,(q-1)\(\frac{\sqrt\pi\;\Gamma(q)}{\Gamma(q+1/2)}\)^{1/q}\,.
\]
%-----------------------------------------------------------------------
\begin{lem}\label{lem:klt}{\rm \cite{MR0121101,Lieb-Thirring76}} Assume that $q\in(1,+\infty)$. Then we have
\[
\Lambda_\R(\mu)=(q-1)^2\,\big(\mu/\mu_1\big)^\beta\quad\forall\,\mu>0\,,
\]
with $\beta=\frac{2\,q}{2\,q-1}$. As a consequence, if $V$ is a nonnegative real valued potential in $\mathrm L^q(\R)$, then we have
\[
\lambda_1^\R[V]\le\Lambda_\R(\nrml Vq)
\]
and equality holds if and only if, up to scalings, translations and multiplications by a positive constant,
\[
V(s)=\frac{q\,(q-1)}{(\cosh s)^2}=:V_1(s)\quad\forall\,s\in\R
\]
where $\nrml{V_1}q=\mu_1$, $\lambda_1^\R[V_1]=\(q-1\)^2$. Moreover the function $\varphi(s)=(\cosh s)^{1-q}$ generates the corresponding eigenspace.
\end{lem}\smallskip
%-----------------------------------------------------------------------
The classical Keller inequality in $\R^d$ asserts that for all $\gamma\ge0$ if $d\ge 3$, $\gamma>0$ if $d=2$, and $\gamma>1/2$ if $d=1$, the lowest negative eigenvalue, $-\lambda_1^{\R^d}[V]$, of the operator $-\,\Delta -V$ satisfies
\[
\lambda_1^{\R^d}[V]^\gamma\le\mathrm L^1_{\gamma,d}\,\nrmRd{V_+}{\gamma+d/2}^{\gamma+d/2}\quad\forall\,V\in\mathrm L^q(\R^d)
\]
where 
\[
\mathrm L^1_{\gamma,d}=\sup\left\{\big(\nrmRd{\nabla u}2^2+\nrmRd u2^2\big)^{-(\gamma+d/2)}:\nrmRd up=1\right\}
\]
is the best constant in the inequality if $q=p/(p-2)$. See~\cite{MR0121101,Lieb-Thirring76,Dolbeault2013437} for details. In the case of infinite cylinders, with $q=\gamma+d/2$ and $\mu=\nrmRd{V_+}a\to+\infty$, the inequality in~$\R^d$ asymptotically determines the \emph{semi-classical regime} for $\Lambda(\mu)$, but another regime appears for cylinders when $\mu>0$ is not too large, as in the case of compact manifolds. This is the content of our main result, Theorem~\ref{Thm:Main}, for which we need one more definition.

Let us introduce
\[
\lambda_\theta:=\(1+\delta\,\theta\,\frac{d-1}{d-2}\)\kappa+\delta\,(1-\theta)\,\lambda_1^{\mathfrak M}\quad\mbox{with}\quad\delta=\frac{n-d}{(d-1)\,(n-1)}\,,
\]
where the dependence on $\theta$ will be discussed at the end of this note, and define
\[
\lambda_\star:=\lambda_{\theta_\star}\quad\mbox{where}\quad\theta_\star:=\frac{(d-2)\,(n-1)\,\big(3\,n+1-d\,(3\,n+5)\big)}{(d+1)\,\big(d\,(n^2-n-4)-n^2+3\,n+2\big)}\,.
\]
%-----------------------------------------------------------------------
\smallskip\begin{thm}\label{Thm:Main} Let $d\ge2$ and $q\in(d/2,+\infty)$. The function $\mu\mapsto\Lambda(\mu)$ is convex, positive and such that
\[\label{asympp}
\Lambda(\mu)^{q-d/2}\sim\mathrm L^1_{q-\frac d2,\,d}\,\mu^q\quad\mbox{as}\quad\mu\to+\infty\,.
\]
Moreover, there exists a positive $\mu_\star$ with
\be{mustar}
\frac{\lambda_\star}{2\,(q-1)}\,\mu_1^\beta\le\mu_\star^\beta\le\frac{\lambda_1^{\mathfrak M}}{2\,q-1}\,\mu_1^\beta
\ee
and $\beta=\frac{2\,q}{2\,q-1}$, such that
\[
\Lambda(\mu)=\Lambda_\R(\mu)\quad\forall\,\mu\in(0,\mu_\star]\quad\mbox{and}\quad\Lambda(\mu)>\Lambda_\R(\mu)\quad\forall\,\mu>\mu_\star\,.
\]
As a special case, if $\M=\S^{d-1}$, inequalities in~\eqref{mustar} are in fact equalities.
\end{thm}\smallskip
%-----------------------------------------------------------------------
In other words, we have shown the Keller-Lieb-Thirring inequality on the cylinder $\mathcal C$:
\[
\tag{KLT}
\lambda_1^{\mathcal C}[V]\le\Lambda\big(\nrmcndgen{V_+}q\big)\quad\forall\,V\in\mathrm L^q(\mathcal C)
\]
where the function $\Lambda:\R^+\to\R^+$ has the properties stated in Theorem~\ref{Thm:Main}. If $\nrmcndgen{V_+}q\le\mu_\star$, optimality is achieved, up to scalings, translations and multiplications by a positive constant, by the potential $V_1$ of Lemma~\ref{lem:klt}. This is based on a rigidity result which, in contrast with results on compact manifolds, involves a non-constant function.

\smallskip
The existence of the function $\mu\mapsto\Lambda(\mu)$ is an easy consequence of a H\"older estimate:
\[
\nrmcndgen{\partial_su}2^2+\nrmcndgen{\nabla_gu}2^2-\iR{V\,|u|^2}\ge\nrmcndgen{\partial_su}2^2+\nrmcndgen{\nabla_gu}2^2-\mu\,\nrmcndgen up^2
\]
with $\mu=\nrmcndgen{V_+}q$ and $q=p/(p-2)$, and of the Gagliardo-Nirenberg-Sobolev inequality
\[
\tag{GNS}
\nrmcndgen{\partial_su}2^2+\nrmcndgen{\nabla_gu}2^2+\alf\,\nrmcndgen u2^2\ge\mu(\alf)\,\nrmcndgen up^2\quad\forall\,u\in\mathrm H^1(\mathcal C)\,,
\]
where $\alf>0$ is a parameter and $\mu(\alf)$ is the corresponding optimal constant. The existence of an optimal function for (GNS) can be argued as in~\cite{Catrina-Wang-01}. Inverting $\alf\mapsto\mu(\alf)$ provides us with $\mu\mapsto\Lambda(\mu)$. See~\cite{DoEsLa} for details and basic properties in a similar case.

The most important point in Theorem~\ref{Thm:Main} is the issue of \emph{symmetry} and \emph{symmetry breaking}. We shall say that there is symmetry if equality in (KLT) is achieved by functions depending only on $s$, and \emph{symmetry breaking} otherwise. By the method used in~\cite{DELT09}, there is a continuous curve $p\mapsto\mu_\star(p)$ defined on $(2,2^*)$ such that symmetry holds if $\mu\le\mu_\star$ and symmetry breaking holds if $\mu>\mu_\star$. It is then clear from the definition of $\Lambda$ and $\Lambda_\star$ that $\Lambda(\mu)\ge\Lambda_\R(\mu)$ with equality if and only $\mu\le\mu_\star$. The main issue is henceforth to estimate $\mu_\star$. Our contribution is based on two perturbation methods:
\begin{itemize}
\item[(i)] For $\mu$ large enough, a non-radial perturbation of an optimal symmetric potential shows symmetry breaking. This is done in the spirit of~\cite{Catrina-Wang-01,Felli-Schneider-03}. The computation gives the upper bound on $\mu_\star$ and a detailed proof is given in Section~\ref{Sec:1d}.
\item[(ii)] For $\mu$ not too large, symmetry holds. A sketch of a proof is given in Section~\ref{Sec:Rigidity}. The key idea is to consider an optimal potential, symmetric or not, and perturb it adequately to prove that it has to be symmetric. The perturbation depends nonlinearly on the minimizer. The proof is not done at the level of (KLT), but at the level of the dual (GNS) inequality. This gives the lower bound on $\mu_\star$. Details will be given in a forthcoming paper,~\cite{DEL2015}.
\end{itemize}

\medskip Apart from Euclidean spaces, very little is known on estimates like the ones of Theorem~\ref{Thm:Main}. A quantitative but non optimal result has been established in~\cite[Corollary 8]{DEL2011}. Some results of symmetry for (KLT) type inequalities have been established for compact manifolds without boundary, see~\cite{DoEsLa,Dolbeault2013437}, and for bounded convex domains in $\R^d$ in relation with the Lin-Ni conjecture, see~\cite{1407}. To our knowledge, the case of non compact manifolds was open so far, apart from the case of the line which was studied in~\cite{Dolbeault06082014} and the partial results of~\cite{DEL2011}. Here we give a result which is optimal when $\M$ is a sphere. Let us finally notice that various observations connecting the sphere, the Euclidean space and the line have been collected in~\cite{Dolbeault06082014}.

%%%%%%%%%%%%%%%%%%%%%%%%%%%%%%%%%%%%%%%%%%%%%%%%%%%%%%%%%%%%%%%%%%%%%%%%
%%%%%%%%%%%%%%%%%%%%%%%%%%%%%%%%%%%%%%%%%%%%%%%%%%%%%%%%%%%%%%%%%%%%%%%%
\section{One-dimensional potentials: proof of Lemma~\ref{lem:klt} and non-symmetric instability}\label{Sec:1d}

We start by a short proof of Lemma~\ref{lem:klt} for the sake of completeness. Notations will be reused in the proof of Lemma~\ref{lem:klt2}.

\medskip\noindent\emph{Proof of Lemma~\ref{lem:klt}.} By applying H\"older's inequality, we get that
\[
\iR{|\partial_su|^2}-\iR{V\,|u|^2}\ge\nrml{\partial_su}2^2-\mu\,\nrml up^2\ge-\,\alf\,\nrml u2^2\,,
\]
where $\mu=\nrml Vq$ and $q=p/(p-2)$. With $V=V_1$, $\alf$ is chosen such that $\mu(\alf)=\mu_1:=\nrml{V_1}q$ where $\mu(\alf)$ is the optimal constant in the inequality
\[
\nrml{\partial_su}2^2+\alf\,\nrml u2^2\ge\mu(\alf)\,\nrml up^2\quad\forall\,u\in\mathrm H^1(\R)\,.
\]
It is standard (see for instance~\cite{Dolbeault06082014}) that the function $\varphi(s)=(\cosh s)^{-2/(p-2)}$ is optimal for the inequality written with $\mu=\mu_1$ and solves
\[
-\,(p-2)^2\,\partial^2_s\,\varphi+\,4\,\varphi-\,2\,p\,\varphi^{p-1}=0\,.
\]
Altogether, this proves Lemma~\ref{lem:klt} when $\mu=\mu_1$ and $\lambda_1[V_1]=\frac4{(p-2)^2}=(q-1)^2$ because $u=\varphi$ and $V=V_1=\frac{2\,p}{(p-2)^2}\,\varphi^{p-2}=q\,(q-1)\,\varphi^{p-2}$ corresponds to the equality case in H\"older's inequality. More details can be found in~\cite{DE2012}.

If $\mu\neq\mu_1$, we can use scalings. Let $V_\nu(s)=\nu^2\,V(\nu\,s)$. If $u_1\neq0$ solves
\[
-\,\partial^2_s\,u_1-\,V\,u_1+\lambda_1[V]\,u_1=0\,,
\]
then $u_\nu(s)=u_1(\nu\,s)$ is an eigenfunction associated with $\lambda_1[V_\nu]=\nu^2\,\lambda_1[V]$. A change of variables shows that $\nrml{V_\nu}q=\nu^{2-1/q}\,\nrml Vq$. Optimality is therefore achieved for~$V_\nu\in\mathrm L^q(\R)$ with $\nrml{V_\nu}q=\mu>0$ if and only if $\nu^{2-1/q}=\mu/\mu_1$ and $V(s)=V_1(s-s_0)$ for some $s_0\in\R$, where
\[
V_{1,\mu}(s)=\nu^2\,V_1(\nu\,s)\quad\forall\,s\in\R\quad\mbox{with}\quad\nu=\big(\mu/\mu_1\big)^\frac q{2\,q-1}\,.
\]
The corresponding eigenfunction is, up to a multiplication by a constant, $\varphi_\mu(s)=\varphi(\nu\,s)$. The lowest eigenvalue for $V\in\mathrm L^q(\R)$ such that $\nrml Vq=\mu$, which realizes the equality in the Keller-Lieb-Thirring inequality, is $\lambda_1[V_{1,\mu}]=\lambda_1[V_1]\,\nu^2=\Lambda_\R(\mu)$. This completes the proof of Lemma~\ref{lem:klt}.\hfill\ \qed\medskip

Next let us consider a function $V$ of $x=(s,z)\in\mathcal C$. Inspired by the results of~\cite{Catrina-Wang-01,Felli-Schneider-03} for Caffarelli-Kohn-Nirenberg inequalities, we can prove that $V_{1,\mu}$, considered as a function on $\mathcal C$, cannot be optimal for the Keller-Lieb-Thirring inequality on $\mathcal C$ if $\mu$ is large enough.
%-----------------------------------------------------------------------
\smallskip\begin{lem}\label{lem:klt2} With the above notations and assumptions, let $V=V(s,z)$ be a nonnegative real valued potential in $\mathrm L^q(\mathcal C)$ for some $q>d/2$ and let $-\,\lambda_1^{\mathcal C}[V]$ be the lowest eigenvalue of the Schr\"odinger operator $-\,\partial^2_s-\,\Delta_g-V$. If $\Lambda_\R(\mu)>4\,\lambda_1^{\mathfrak M}/(p^2-4)$, then
\[
\sup\left\{\lambda_1^{\mathcal C}[V]\,:\,V\in\mathrm L^q(\mathcal C)\,,\;\nrmcndgen Vq=\mu\right\}>\Lambda_\R(\mu)\,,
\]
that is, the above maximization problem cannot be achieved by a potential $V$ depending only on the variable~$s$.
The condition $\Lambda_\R(\mu)>4\,\lambda_1^{\mathfrak M}/(p^2-4)$ is explicit and equivalent to $\mu^\beta>\mu_1^\beta\,\lambda_1^{\mathfrak M}/(2\,q-1)$.
\end{lem}\smallskip
%-----------------------------------------------------------------------
\begin{proof} Let $\varphi_\mu$ be as in the proof of Lemma~\ref{lem:klt}. We argue by contradiction and consider
\[
\phi_\varepsilon(s,z):=\varphi_\mu(s)\,+\varepsilon\,\big(\varphi_\mu(s)\big)^{p/2}\,\psi_1(z)\quad\mbox{and}\quad V_\varepsilon(s,z):=\mu\,\frac{|\phi_\varepsilon(s,z)|^{p-2}}{\nrmcndgen{\phi_\varepsilon}p^{p-2}}\,.
\]
Here $\psi_1$ denotes an eigenfunction of $-\,\Delta_g$ on $\M$ such that $\nrM{\psi_1}2=1$ and $-\,\Delta_g\psi_1=\lambda_1^{\mathfrak M}\,\psi_1$. Then
\begin{multline*}
\nrmcndgen{\partial_s\phi_\varepsilon}2^2+\nrmcndgen{\nabla_g\phi_\varepsilon}2^2-\,\mu\,\nrmcndgen{\phi_\varepsilon}p^2-\(\nrml{\partial_s\varphi_\mu}2^2-\,\mu\,\nrml{\varphi_\mu}p^2\)\\
=\varepsilon^2\(\lambda_1^{\mathfrak M}\,\nrml{\varphi_\mu}p^p+\nrml{\partial_s\varphi_\mu^{p/2}}2^2-\,(p-1)\,\mu\,\nrml{\varphi_\mu}p^{2-p}\iR{\varphi_\mu^{2\,(p-1)}}\)+o(\varepsilon^2)\,.
\end{multline*}
Since $\mu\,\nrml{\varphi_\mu}p^{2-p}\,\varphi_\mu^{p-2}=V_{1,\mu}(s)=\nu^2\,V_1(\nu\,s)$ with $\nu=\big(\mu/\mu_1\big)^\beta$ as in the proof of Lemma~\ref{lem:klt}, since $\chi=\varphi^{p/2}$ solves
\[
-\,\partial^2_s\,\chi-\,(p-1)\,V_1\,\chi=-\,\Big(\frac p{p-2}\Big)^2\chi\,,
\]
and since $(p-1)\,q\,(q-1)=\frac{2\,p\,(p-1)}{(p-2)^2}$, we get that
\[
\nrml{\partial_s\varphi_\mu^{p/2}}2^2-\,\frac{2\,p\,(p-1)}{(p-2)^2}\,\nu^2\iR{\varphi_\mu^{2\,(p-1)}}=-\,\nu^2\,\Big(\frac p{p-2}\Big)^2\,\nrml{\varphi_\mu}p^p\,.
\]
We recall that $\nu^2=\Lambda_\R(\mu)/(q-1)^2=\frac14\,(p-2)^2\,\Lambda_\R(\mu)$. According to~\cite[Appendix~A]{Dolbeault06082014}, we see that
\[
\nrml{\varphi_\mu}p^p=\frac4{p+2}\,\nrml{\varphi_\mu}2^2\quad\mbox{and}\quad\nrmcndgen{\phi_\varepsilon}2^2=\Big(1+\frac{4\,\varepsilon^2}{p+2}\Big)\,\nrml{\varphi_\mu}2^2\,.
\]
Hence we finally find that
\[
-\,\lambda_1^{\mathcal C}[V_\varepsilon]+\,\Lambda_\R(\mu)\le\frac{4\,\varepsilon^2}{p+2}\(\lambda_1^{\mathfrak M}-\,\tfrac14\,(p^2-4)\,\Lambda_\R(\mu)\)+o(\varepsilon^2)\,.
\]
This shows that $\lambda_1[V_{1,\mu}]-\lambda_1^{\mathcal C}[V_\varepsilon]<0$ is negative for $\varepsilon>0$, small enough, if \hbox{$\Lambda_\R(\mu)>\frac{4\,\lambda_1^{\mathfrak M}}{p^2-4}$}.\end{proof}

The condition found by V.~Felli and M.~Schneider in~\cite{Felli-Schneider-03} can be recovered by noticing that $\lambda_1^{\S^{d-1}}=d-1$, when $\M=\S^{d-1}$. In that case, the above computation are exactly equivalent to the computations for Caffarelli-Kohn-Nirenberg inequalities: see~\cite{Dolbeault06082014} for details.

%%%%%%%%%%%%%%%%%%%%%%%%%%%%%%%%%%%%%%%%%%%%%%%%%%%%%%%%%%%%%%%%%%%%%%%%
%%%%%%%%%%%%%%%%%%%%%%%%%%%%%%%%%%%%%%%%%%%%%%%%%%%%%%%%%%%%%%%%%%%%%%%%
\section{Symmetry: a rigidity result}\label{Sec:Rigidity}

In this section we get a lower bound on $\mu_\star$ and complete the proof of Theorem~\ref{Thm:Main}. Let us define
\[
\mathsf J[V]:=\frac{\nrmcndgen Vq^q-\nrmcndgen{\partial_sV^{(q-1)/2}}2^2-\nrmcndgen{\nabla_gV^{(q-1)/2}}2^2}{\nrmcndgen{V^{(q-1)/2}}2^2}\,.
\]
We shall consider a critical point of $\mathsf J$ and prove that it is symmetric using a well chosen perturbation.
%-----------------------------------------------------------------------
\smallskip\begin{lem}\label{lem:equiv} With the notations of Section~\ref{Sec:Intro} and under the assumptions of Theorem~\ref{Thm:Main}, we have
\[
\Lambda(\mu)=\sup\left\{\mathsf J[V]:\nrmcndgen Vq=\mu\right\}\,.
\]\end{lem}\smallskip
%-----------------------------------------------------------------------
\noindent\emph{Sketch of the proof of Lemma~\ref{lem:equiv}.} There exists a nonnegative potential $V\in\mathrm L^q\cap C^\infty(\mathcal C)$ such that $\mathsf J[V]=\Lambda(\mu)$. Lemma~\ref{lem:equiv} is based on the equivalence of (KLT) and (GNS), which can also be seen by using $u=V^{(q-1)/2}$ and considering the equality case in H\"older's inequality. See~\cite{DEL2015} for more details.
\hfill\ \qed\medskip

To any potential $V\ge0$ we associate the \emph{pressure} function
\[
\mathsf p_V(r):=r\,V(s)^{-\frac12}\quad\forall\,r=e^{-\alpha s}\,.
\]
With $\alpha=\frac1{q-1}\,\sqrt{\Lambda_\R(\mu)}$, let us define
\begin{multline*}
\mathsf K[\mathsf p]:=\frac{2q-1}{2\,q}\,\alpha^4\irdmu{\left|\mathsf p''-\frac{\mathsf p'}r-\frac{\Delta_g\mathsf p}{\alpha^2\,(2q-1)\,r^2}\right|^2\mathsf p^{1-2q}}\\
+2\,\alpha^2\, \irdmu{\frac1{r^2}\,\left|\nabla_g\mathsf p'-\frac{\nabla_g \mathsf p}r\right|^2\mathsf p^{1-2q}}+\(\lambda_\star-\frac2{q-1}\,\Lambda_\R(\mu)\)\irdmu{\frac{|\nabla_g\mathsf p|^2}{r^4}\,\mathsf p^{1-2q}}\,,
\end{multline*}
where $d\mu$ is the measure on $\R^+\times\M$ with density $r^{2q-1}$, and $'$ denotes the derivative with respect to $r$.%-----------------------------------------------------------------------
\smallskip\begin{lem}\label{lem:K} There exists a positive constant $\mathsf c$ such that, with the above notations, if $V$ is a critical point of $\mathsf J$ under the constraint $\nrmcndgen Vq=\mu$, then $\mathsf K[\mathsf p_V]=0$.\end{lem}\smallskip
%-----------------------------------------------------------------------
The sketch of the proof is as follows. Because $V$ is critical point, we know that
\[
0=\mathsf J'[V]\cdot V^{1-q}\,\left[\partial_s\(e^{-2(q-1)\alpha s}\,\partial_s\(e^{(2q-1)\alpha s}\,V^\frac{2q-1}2\)\)+e^{\alpha s}\,\Delta_g V^\frac{2q-1}2\right]\,.
\]
On the other hand, it has been shown in~\cite{DEL2015} that the r.h.s.~is greater or equal than a positive constant times $\mathsf K[\mathsf p_V]$.

\medskip\noindent\emph{Sketch of the proof of Theorem~\ref{Thm:Main}.} If $\mu<\mu_\star$, then $\lambda_\star-\frac2{q-1}\,\Lambda_\R(\mu)$ is positive and we read from the expression of $\mathsf K[\mathsf p_V]$ that $\mathsf p_V$ depends only on $r$ or, equivalently, $V$ depends only on $s$. Up to scalings, translations and multiplications by a positive constant, we get that $V=V_1$ as defined in Lemma~\ref{lem:klt}.
\hfill\ \qed\medskip

Note that the constant $\lambda_\star$ is an estimate of the largest constant $\lambda$ such that
\[
\iM {\(\tfrac12\,\Delta_g(|\nabla_g f|^2)-\nabla_g(\Delta_g f)\cdot\nabla_gf-\tfrac1{2q-1}\,(\Delta_g f)^2-\lambda\,|\nabla_g f|^2\)\,f^{1-2q}}\ge0\,,
\]
for any positive function $f\in C^3(\M)$. It is estimated by $\lambda_\theta$ with $\theta\in[\theta_\star,1]$, according to~\cite{DEL2015}. In the case of the sphere, that is, $\mathfrak M=\S^{d-1}$, we have that $\frac{d-1}{d-2}\,\kappa=\lambda_1^{\mathfrak M}$ and $\lambda_\theta=\big(1+\delta\,\frac{d-1}{d-2}\big)\,\kappa=\big(\frac{d-2}{d-1}+\delta\big)\,\lambda_1^{\mathfrak M}$ is independent of $\theta$. Otherwise, by Lichnerowicz' theorem, we know that $\frac{d-1}{d-2}\,\kappa\le\lambda_1^{\mathfrak M}$. Hence $\theta\mapsto\lambda_\theta$ is a non-increasing function, and since~$\theta_\star$ is always negative, we have a simple lower bound for $\lambda_\star$:
\[
\lambda_\star\ge\lambda_0=\kappa+\delta\,\lambda_1^{\mathfrak M}\,.
\]

%%%%%%%%%%%%%%%%%%%%%%%%%%%%%%%%%%%%%%%%%%%%%%%%%%%%%%%%%%%%%%%%%%%%%%%%
%%%%%%%%%%%%%%%%%%%%%%%%%%%%%%%%%%%%%%%%%%%%%%%%%%%%%%%%%%%%%%%%%%%%%%%%

\smallskip{\small\noindent{\bf Acknowledgments.} J.D.~has been partially supported by the projects \emph{STAB} and \emph{Kibord} (J.D.) of the French National Research Agency (ANR). M.L.~has been partially supported by the US - NSF grant DMS-1301555.\\[2pt]
{\sl\small \copyright~2015 by the authors. This paper may be reproduced, in its entirety, for non-commercial purposes.}
\begin{flushright}{\sl\today}\end{flushright}}

%%%%%%%%%%%%%%%%%%%%%%%%%%%%%%%%%%%%%%%%%%%%%%%%%%%%%%%%%%%%%%%%%%%%%%%%
%%%%%%%%%%%%%%%%%%%%%%%%%%%%%%%%%%%%%%%%%%%%%%%%%%%%%%%%%%%%%%%%%%%%%%%%
%\newpage
%\bibliographystyle{abbrv}\bibliography{../References}

%%%%%%%%%%%%%%%%%%%%%%%%%%%%%%%%%%%%%%%%%%%%%%%%%%%%%%%%%%%%%%%%%%%%%%%%

%%%%%%%%%%%%%%%%%%%%%%%%%%%%%%%%%%%%%%%%%%%%%%%%%%%%%%%%%%%%%%%%%%%%%%%%
\end{document}